\documentclass[]{article}
\usepackage[]{amsmath, amsthm, amsfonts}

\newcommand {\N} {{\mathbb N}}
\newcommand {\C} {{\mathbb C}}

\newcommand {\Z} {{\mathbb Z}}
\newcommand {\Q} {{\mathbb Q}}
\newcommand {\PP} {{\mathbb P}}
\newcommand {\HH} {{\mathcal H}}
\newcommand {\qk} {quasi-K\"ahler}

\newtheorem{thm}[subsection]{Theorem}
\newtheorem{cor}[subsection]{Corollary}
\newtheorem{lemma}[subsection]{Lemma}
\newtheorem{prop}[subsection]{Proposition}
\newtheorem{defn}[subsection]{Definition}
\newtheorem{remark}[subsection]{Remark}

\title{
	Solvable fundamental groups of algebraic varieties and K\"ahler 
	manifolds
}
\author{
	Donu Arapura \thanks{Partially  supported by NSF}\\
	Department of Mathematics\\ Purdue University \\ West Lafayette, IN 
	47907
	 \and
	Madhav Nori\\
	Department  of Mathematics \\ University of Chicago\\ Chicago, IL 
	60637
}

\date{ June 24, 1997}

\begin{document}
\maketitle

Our objective, in this paper, is to gain some  understanding  of 
those groups which arise as fundamental groups
of  compact K\"ahler manifolds, Zariski open subsets of  compact
K\"ahler manifolds, or normal complex algebraic varieties.
The groups of the first type have come to be known as  K\"ahler
groups, and we will refer to those of the second as {\qk}.
While the complete structure of these  groups seems rather 
intractable at the moment
(see \cite{KahlerBook} for the state of of the art), 
the structure of certain subclasses are becoming much clearer.
Of specific interest for us is the class of nilpotent groups.
Many restrictions on, as well as interesting examples of, nilpotent  K\"ahler groups 
have been found  by Campana \cite{C} and Carlson and Toledo \cite{CT}.
 Although nilpotent {\qk} groups have not been studied 
 systematically, nontrivial constraints can be obtained from Morgan's work \cite{Morgan}.
The main conclusion of this paper, is that if one casts the net a 
little wider, then no really interesting new examples are obtained.
In particular, we will show that a  polycyclic {\qk} group is
 virtually nilpotent, which is to say that it contains a nilpotent subgroup of
 finite index. For algebraic varieties, we
can make an even stronger statement that  the fundamental 
group of a normal variety must be virtually nilpotent if
it is  solvable and possesses a faithful representation into $GL_n(\Q)$

The first two sections of this paper are purely group theoretic. 
We introduce the class of solvable groups of finite rank, which
contains the class of   solvable
subgroups of $GL_n(\Q)$ (and polycyclic groups in particular).
 Canonically attached to every  such group
$\Gamma$ is an algebraic group $H(\Gamma)$ defined over $\Q$,
and a representation $\Gamma\to H(\Gamma)$ with Zariski dense
image.
This generalizes earlier constructions of Malcev \cite{Malcev} and
Mostow \cite{M} for nilpotent and polycyclic groups respectively.
The third section contains the first main theorem, that if 
the quotient $\Gamma$ of the   fundamental group of a normal
variety by a term of the derived series has finite rank, then
$\Gamma$ contains an extension of a nilpotent group by a torsion
group as a subgroup of finite index. The basic idea is to choose
a finitely generated field of definition, and
observe that the Galois group acts on the group of $\Q_l$ points of $H(\Gamma)$
for almost all primes $l$. This together with certain arithmetic 
considerations, forces the identity component $H(\Gamma)^o$ to be
unipotent, and the theorem follows easily from this. The
fourth section of the paper contains the second main theorem:
If $\pi$ is a {\qk } group,  with finitely generated derived group,
such that the quotient $\Gamma$ of $\pi$ by a term of the
derived series has finite rank, then $\Gamma$ contains a nilpotent
by torsion subgroup of finite index. Once again the strategy is
to establish the 
unipotency of $H(\Gamma)^o$. However,  this time it  is reduced  to a 
homological statement  which is shown to be a consequence
a generalization of a theorem
of Beauville \cite{B} obtained by the first author \cite{A2}.

The reader will certainly have noticed that the two theorems
are very similar in content, but  quite different in
methods of proof.  This by itself should not come as a surprise; 
a number of results in algebraic geometry, such as the
quasiunipotence of local monodromy, have been proved
by both arithmetic and transcendental methods.
In fact, the style
of the first argument is quite similar to that of Grothendieck's proof this
 theorem. What does seem a bit curious is the slight disparity of the results.
In particular,  there exists  solvable
subgroups of $GL_n(\Q)$ which violate the conclusion of the first theorem,
but not the second. 
We leave open the question of whether
such groups can be \qk. 
It is worth noting that a  second proof that  nonvirtually nilpotent solvable
subgroups of $GL_n(\Q)$ are not fundamental groups of smooth 
projective varieties can be obtained by combining the above arguments
with those of Simpson \cite{Si}. Once again, Simpson's arguments are
arithmetic in nature and do not apply to nonalgebraic K\"ahler manifolds.

This paper has existed, for a long time, in a chimerical state:
occasionally referenced but never seen. We apologize for that.
We would like to thank F. Campana, J. Carlson, and D. Toledo for sending us 
their then preprints, and M. S. Raghunathan for pointing out the reference 
to \cite{M}.

\section{Preliminaries on Algebraic Groups}

A general reference for this section is \cite{Bo}. Let $DG = D^1G$ be the 
derived subgroup of a group $G$, and set $D^iG= DD^{i-1}G$. 
If $F \subset E$ is an extension of fields, and $G$ is an algebraic 
group defined over $F$, let $G_{E} = G\times_{spec\, F}spec\, E$. 
For the remainder of this section, $G$ will denote an algebraic
group over a field $F$ of characteristic $0$, and $U(G)$ will be its
unipotent radical. 
We will make use of the following result of Mostow \cite{M56}:

\begin{thm}\label{thm:splitting}
If $G$ is as  above, then the exact sequence:
$$1\to U(G)\to G\to G/U(G)\to 1$$
is split, and any two splittings are conjugate by an $F$ rational point
of $U(G)$.
\end{thm}

Let $V$ be the centralizer of $U(G)$ in $G$, and let $W= V\cap U(G)$.
Clearly $W$ is the unipotent radical of $V$. It follows that 
$$1\to W\to V\to V/W\to 1$$
has a unique splitting $s: V/W\to V$, because $W$ central in $V$.
Put $N(G) = s(V/W)$. By construction $N(G)$ is reductive and invariant 
under all automorphisms of $G$; in particular, it is normal in $G$.

\begin{lemma}\cite[lemma 4.6]{M} Every normal reductive
algebraic subgroup of $G$ is contained in $N(G)$.
\end{lemma}

\begin{defn} An algebraic group $G$ is minimal if $N(G)$ is trivial.
Let $G_{min} = G/N(G)$; then this is a minimal algebraic group. An 
algebraic group $G$ is minimal if and only if the centralizer of 
$U(G)$ is contained in $U(G)$, as we see from the definition of $N(G)$.
\end{defn}

\begin{lemma}\label{lemma:a3}
 Let $f:G\to G'$ be a homomorphism of algebraic groups 
for which $f(G)$ is normal in $G'$. Then $f(N(G)) \subseteq N(G')$,
and consequently, $f$ induces a homomorphism $f_{min}:G_{min}\to G_{min}'$.
If $f$ is a injection (respectively surjection), then $f_{min}$
is also a injection (respectively surjection).
\end{lemma}

\begin{proof}
As $f(N(G))$ is normal and reductive in $f(G)$, we see that 
$f(N(G))\subseteq N(f(G))$. As $N(f(G))$ is invariant under all 
automorphisms of $f(G)$, and in particular the restriction of inner 
automorphisms of $G'$ to $f(G)$, it follows that $N(f(G))\subseteq 
N(G')$. This proves the first assertion.

Because $H = f(G)\cap N(G')$ is normal in the reductive group $N(G')$,
we see that $H$ is reductive. Also, $H$ is normal in $f(G)$, and this 
implies that $H=N(f(G))$, Thus $ker(f) = \{1\}$ implies $\ker(f_{min})
=\{1\}$. 
\end{proof}

\begin{lemma}\label{lemma:a4} Let $G$ be a minimal algebraic group defined over a 
field $F$ of characteristic $0$. There is an affine algebraic group 
$A$, defined over $F$, 
that acts on $G$ so that for all fields $E\supseteq F$, $A(E)\to 
Aut(G_E)$ is an isomorphism.
\end{lemma}

\begin{proof} If $G = U(G)$, the automorphism of $G$ are just 
automorphisms of its Lie algebra, and these evidently form an affine
algebraic group, denoted by $Aut(G)$.

In the general case, by \ref{thm:splitting} we can
express $G$ as a semidirect product of $M$ and 
$U(G)$, where $M$ is reductive. The homomorphism $\rho: M\to Aut(U(G))$
is faithful, because $G$ is minimal. Let $X$ be the normalizer of 
$\rho(M)$ in $Aut\, U(G)$. Thus $X$ is an algebraic group, and $X$ acts 
naturally on $U(G)$ and $M$, and hence on $G$, their semidirect product.
Now $G$ acts on itself by conjugation. Thus $Y$, the semidirect 
product of $X$ and $G$, acts on $G$. For any field $E\supset F$,
we see $X(E)\to Aut G_E$ is one to one and its image equals
$\{\phi\in Aut G_E\, |\, \phi(M_E) = M_E\}$. But, if $\phi\in 
Aut(G_E)$, $\phi(M_E)$ is a conjugate of $M_E$, and we deduce that 
$Y(E)\to Aut G_E$ is surjective for all fields $E$ containing $F$.
Finally, the coordinate ring $R$ of $G$ is generated as an $F$-algebra
by a finite dimensional $Y$-stable subspace $V\subset R$. Put $K = 
ker(Y\to GL(V))$ and let $A = Y/K$. We see that $A$ is the desired 
algebraic group.

\end{proof}

The algebraic group $A$ in the previous lemma will be denoted by $Aut\, 
G$ henceforth.

\begin{lemma}\label{lemma:a5} If $G$ is a solvable minimal algebraic group, 
then the  action of $Aut(G)^o$ on $G^o/U(G)$ is trivial.

\end{lemma}

\begin{proof} This follows immediately from the fact that the action 
of a connected algebraic group on a torus is trivial. 
\end{proof}

\begin{lemma}\label{lemma:a6} Let $\Lambda$ be  a directed set, and 
$\{G_\lambda\}_{\lambda\in \Lambda}$ 
a directed system of minimal algebraic groups, such that for each $\lambda\le 
\mu$, we have $G_\lambda \subseteq G_\mu$  and $U(G_\lambda) =
U(G_\mu)$. Then there is a minimal algebraic group $G$ and a 
monomorphism $f_\lambda: G_\lambda \to G$ for each $\lambda\in 
\Lambda$, so that $f_\lambda = f_\mu |_{G_\lambda}$ whenever 
$\lambda \le \mu$.
\end{lemma}

\begin{proof}
 Let $U= U(G_\lambda)$ for all $\lambda\in \Lambda$.
Let 

$$S_\lambda = \{M\subseteq G_\lambda\, |\, 
M\> is\> a\> closed\> subgroup\> and\>
M\to G_\lambda/U \> is \> an \> isomorphism\}$$

Then $U$ acts transitively on $S_\lambda$. Also $M\to M\cap G_\lambda$
gives a $U$-equivariant morphism from $S_\mu$ to $S_\lambda$ whenever
$\lambda \le \mu$. Now choose $\lambda_0\in \Lambda$ so that $dim\, 
S_{\lambda_0} \ge dim\, S_\lambda$ for all $\lambda \in \Lambda$.
From the above, we see that $S_\mu \to S_{\mu_0}$ is a bijection if
$\lambda_0 \le \mu$. It follows that there is a collection
$\{M_\lambda\, | \, \lambda \in \Lambda \}$ with $M_\lambda \in 
S_\lambda$, and $M_\mu \cap G_\lambda = M_\lambda$ whenever $\lambda 
\le \mu$.

If $\rho_\lambda : M_\lambda \to Aut U$ denotes the conjugation action
of $M_\lambda$ on $U$, we have seen that $\rho_\lambda$ is one to one
because $G_\lambda$ is minimal. Also, the inequality $\lambda \le \mu$
implies that $\rho_\lambda (M_\lambda) \subseteq \rho_\mu (M_\mu)$.
Let $M$ be the Zariski closure of $\cup \{\rho_\lambda(M_\lambda)\, |\,
\lambda \in \Lambda \} $. This is reductive. Let $G$ be the semidirect
product of $M$ and $U$, and define $f_\lambda: G_\lambda \to M$ by
$f_\lambda(u) = u$ for all $u\in U$, and  $f_\lambda = \rho_\lambda$
on $M_\lambda$. This completes the proof.

 \end{proof}


Let $G$ be a connected solvable group defined over $F$. 
Let $U=U(G)$  and  $T$ a maximal torus.
As noted previously, $Aut\,U$ is an algebraic group which coincides
with the group of automorphisms of the Lie algebra $N$ of $U$.
There is a homomorphism of algebraic groups $T\to Aut\, U$ given
by conjugation. $G$ is the semidirect product of $T$ with $U$.

 \begin{lemma} With the previous notation,
 suppose that $T$ acts trivially (by conjugation) on $U/DU$
 then $G$ is isomorphic to the product of $U$ and $T$, and is 
 therefore nilpotent.
 \end{lemma}
 
 \begin{proof} We have to show that 
 $T$ acts trivially on $U$.  
 Let $S$ be the subgroup of $Aut(N)$ of automorphisms $\sigma$ satisfying
 $(1-\sigma)(N) \subseteq [N,N]$. The elements of $S$ are unipotent.
 By assumption, the image of the homomorphism $T \to Aut(U)= Aut(N)$ lies in 
 $S$ and so the map must be trivial.
 \end{proof}

There is a homomorphism of algebraic groups $G/DG\to Aut\, DG/D^{2}G$
given by conjugation.

\begin{lemma}\label{lemma:unip}
 With the previous notation,
 suppose that image of $G/DG$ in $Aut\linebreak[0] DG/D^2G$ is unipotent.
  Then $G$ is the product of  $U$ with  $T$.
\end{lemma}
  
 \begin{proof} Note that $G/DG$ is a product of a reductive group, 
 which is isomorphic to $T$, and a unipotent group $U'$, which is 
 isomorphic to the image of $U$ under the projection $G\to G/DG$.
 The action of $T$ on $U'$ by conjugation is of course trivial.
 Consider the exact sequence
$$ DG/D^2G \to U/DU \to U' \to 0.$$
By assumption, the image of $T\subset G/DG$ in $Aut\, DG/D^2G$ is 
unipotent, and therefore trivial. 
Thus $T$ acts trivially on the image of $DG/D^2G$ in $U/DU$ as
well as on it the cokernel.  
Therefore $T$  acts trivially on $U/DU$, and the lemma follows from
the previous one.
\end{proof}

\section{ Solvable groups of finite rank}

In this section, we shall associate to a solvable group $\Gamma$ of
finite rank, an algebraic group of $H(\Gamma)$ defined over $\Q$
and a homomorphism $i(\Gamma):\Gamma\to H(\Gamma)(\Q)$ with Zariski
dense image and a torsion subgroup as kernel. While $\Gamma \mapsto 
H(\Gamma)$ is a not a functor, automorphisms of 
$\Gamma$ will extend to automorphisms of $H(\Gamma)$.  Furthermore when
$\Gamma$ is finitely generated, automorphisms of its profinite 
completion $\hat \Gamma$ extend to automorphisms of $H(\Gamma)_{\Q_l}$
for all but finitely many primes $l$. 
Our construction of $H(\Gamma)$ is identical to Mostow's \cite{M}
in the case where $\Gamma$ is polycyclic, but our
proof that this construction works for solvable groups of finite rank
is necessarily a bit more complicated. Also the results on the 
profinite completion are used crucially when Galois theory is applied.
For these reasons, we have chosen to give all the details of the proofs.

Let us recall the standard construction of the proalgebraic hull of
$\HH(\pi, F)$ associated to a topological group $\pi$ and a 
topological field $F$. One considers the category $C(\pi, F)$ where the
objects are pairs $(G,f)$ with $G$ an affine algebraic group defined 
over $F$, and $f:\pi \to G(F)$ a continuous homomorphism with Zariski 
dense image. A morphism $(G,f)\to (G', f')$ in our category is simply 
a commutative diagram
$$\begin{array}{ccccc}
G&         &\stackrel{\phi}{\longrightarrow}&        & G'\\
 &\nwarrow &                               & \nearrow& \\
 &         &             \pi               &         &     \\
\end{array}
$$
such that $\phi$ is homomorphism of algebraic groups. Such a 
$\phi$, if it exists is unique and necessarily an epimorphism because
$f$ and $f'$ have Zariski dense images. In particular, by
lemma \ref{lemma:a3}, $\phi_{min}: G_{min} \to G_{min}'$ is defined.
Set
$$\HH(\pi, F) = \lim_{\stackrel{\longleftarrow}{(G,f)}}\, G$$
and 
$$H(\pi, F) = \lim_{\stackrel{\longleftarrow}{(G,f)}}\, G_{min}.$$

Some easy observations follow.

\begin{remark}\label{remark:1} $\pi \mapsto \HH(\pi, F)$ is a functor, but $\pi \mapsto
H(\pi, F)$ is not. However, if $a:\pi\to \pi'$ is a continuous 
homomorphism with the closure of $a(\pi)$ normal in $\pi'$, we deduce,
from lemma \ref{lemma:a3} that $H(a,F):H(\pi, F)\to H(\pi',F)$ is
defined. In particular, if $\Gamma$ is a discrete group and $t: 
\Gamma\to \hat\Gamma$ is the homomorphism to its profinite completion,
we have a natural epimorphism 
$$H(t,F): H(\Gamma, F) \to H(\hat \Gamma, F)$$
\end{remark}

\begin{remark} If $F\to E$ is a continuous homomorphism of fields, we 
have a functor $C(\pi, F) \to C(\pi, E)$. And this induces an 
epimorphism $H(\pi,E)\to H(\pi, F)_E$.
\end{remark}

\begin{remark}\label{remark:3} By construction, we have a continuous homomorphism
$i(\pi, F)$ from $\pi$ to the group of $F$-rational points of
$H(\pi,F)$. With $a:\pi\to\pi'$ as in remark \ref{remark:1}, we have:
$$H(a,F)\circ i(\pi, F) = i(\pi', F)\circ a$$
\end{remark}

\begin{defn} A solvable group $\Gamma$ has finite rank,
if there is a decreasing sequence
$$\Gamma = \Gamma_0\supset \Gamma_1 \supset \ldots \supset 
\Gamma_{m+1} = \{1\}$$
of subgroups, each normal in its predecessor,
 such that $\Gamma_i/\Gamma_{i+1}$ is abelian and
$\Q \otimes (\Gamma_i/\Gamma_{i+1})$ is finite dimensional for all
$i$. The rank 
$$rk(\Gamma) = \sum_{i=0}^m\, rk(\Q \otimes 
(\Gamma_i/\Gamma_{i+1}))$$
is clearly independent of the choice of the sequence $\{\Gamma_i\}$.
\end{defn}

This is a weakening of the notion of a polycyclic group which,
in the above terms, amounts to requiring that each $\Gamma_i/\Gamma_{i+1}$
is finitely generated.

For the remainder of this section, all groups considered are solvable
of finite rank with discrete topology unless indicated otherwise.
 We endow $\Q$ with discrete
topology, and abbreviate $H(\pi, \Q), i(\pi,\Q), H(a, \Q)$ by
$H(\pi)$ etcetera. The only fields $F$ considered have characteristic
zero.

\begin{thm}\label{thm:solv} Let $\Gamma$ be a solvable group of finite rank with 
discrete topology. Then

(A) $H(\Gamma)$ is an algebraic group (and not just a proalgebraic group).

(B) $rk(\Gamma) = dim \, U(H(\Gamma))$.

(C) The kernel of $i(\Gamma) : \Gamma \to H(\Gamma)(\Q)$ is a torsion
group.

(D) The image of $i(\Gamma)$ is Zariski dense.

\end{thm}

Part (D) of the theorem is immediate from the construction. When 
$\Gamma$ is polycyclic, the theorem is due to Mostow \cite[4.9]{M}.
As a corollary we obtain a natural characterization of these groups.

\begin{cor} A solvable group $\Gamma$ has finite rank if and only if there 
exists a torsion normal subgroup $N\subset\Gamma$ such that 
$\Gamma/N$ possesses an embedding into an affine algebraic group
defined over $\Q$.
\end{cor}

\begin{proof} One direction follows from the theorem. For the 
converse, suppose that $N\subseteq \Gamma$ is torsion subgroup
and $i:\Gamma\to G(\Q)$ a homomorphism into an algebraic
group with kernel $N$. We can assume that $G$ is solvable, 
after replacing it by the Zariski closure of $\Gamma$. Then the
sequence $\Gamma_i= i^{-1}(D^iG(\Q))$ has the required properties.
\end{proof}

\begin{defn} $n(\Gamma, F) = 
sup\,\{ dim\, U(G) \,|\, (G,f) \in Obj\, C(\Gamma, F) \}.$
Thus $n(\Gamma, F) \in \N \cup \{\infty\}$.
\end{defn}

\begin{lemma} \label{lemma:s6}
If $1\to \Gamma' \to \Gamma \to \Gamma'' \to 1$ is 
exact, $n(\Gamma, F) \le n(\Gamma', F) + n(\Gamma'', F)$.
\end{lemma}

\begin{proof} Let $(G,f)$ be an object of $C(\Gamma, F)$. Now let
$G' $ be the Zariski closure of $f(\Gamma')$ and let $G'' = G/G'$.
Then $f$ induces  $\bar f:\Gamma'' \to G''(F)$ and both 
$(G',f|_{\Gamma'})$ and $(G'', \bar f)$ are in $C(\Gamma',F)$
and $C(\Gamma'', F)$ respectively.  A short exact sequence of algebraic
groups induces a short exact sequence of unipotent radicals, so the 
lemma follows.
\end{proof}

\begin{lemma}\label{lemma:s7} $n(\Gamma, F) \le rk(\Gamma)$, where $F$
is a field of characteristic zero.
\end{lemma}

\begin{proof}

From the previous lemma, by induction on the length of the derived 
series of $\Gamma$, we are reduced to the case where $\Gamma$ is 
abelian. If $(G,f) \in Obj\, C(\Gamma, F)$, then the inclusion
$U(G) \hookrightarrow G$ is split by a homomorphism
$p:G \to U(G)$. 
But in this case, $U(G)$ is a vector space, spanned by $p(f(\Gamma))$,
and so evidently $rk(\Q\otimes \Gamma) \ge dim\, U(G)$.

\end{proof}

\begin{lemma}\label{lemma:s8} If $(G, f)\in Obj\, C(\Gamma, F)$ satisfies
$n(\Gamma, F) = dim\, U(G)$, then $H(\Gamma, F) \to G_{min}$
is an isomorphism. In particular, $H(\Gamma, F) $ is an 
algebraic group and $dim\, U(H(\Gamma, F)) = n(\Gamma, F)$.
\end{lemma}

\begin{proof} For the first assertion, we need to check that any
morphism $\phi: (G', f') \to (G, f)$ induces an isomorphism
$\phi_{min}: G_{min}' \to G_{min}$. Now, $\phi$ being as epimorphism,
restricts to an epimorphism of unipotent radicals. This gives
$$n(\Gamma, F) \ge dim\, U(G') \ge dim \, U(G) = n(\Gamma, F).$$

Therefore, $U(G') \to U(G)$ gives an isomorphism of Lie algebras, and
is itself an isomorphism. Consequently, $ker\phi$ is reductive, and
$\phi_{min} $ is an isomorphism.

For the second assertion, we need to know that such a $(G,f)$ exists,
and this is assured by the previous lemma.
\end{proof}

\begin{lemma}\label{lemma:s9} Let $1\to \Gamma'\to \Gamma \to \Gamma'' \to 1$ be exact.
Assume that $i(\Gamma', F)$ extends to a homomorphism $j:\Gamma\to 
G(F)$ where $G$ is an algebraic group that contains $H(\Gamma', F)$.
Then

(A) $H(\Gamma', F)\to H(\Gamma, F)$ is a monomorphism.

(B) $n(\Gamma, F) = n(\Gamma', F) + n(\Gamma'', F)$

(C) $1\to ker\, i(\Gamma', F)\to 
ker\, i(\Gamma, F)\to ker\, i(\Gamma'', F)$
is exact.

\end{lemma}

\begin{proof} Let $q:\Gamma' \to \Gamma$ and $p:\Gamma' \to \Gamma''$ 
denote the given homomorphisms. By remark \ref{remark:1}, we have
$H(q,F):H(\Gamma',F)\to H(\Gamma, F)$ and 
$H(p,F):H(\Gamma,F)\to H(\Gamma'', F)$. Now $H(p,F)$ is an epimorphism
whose kernel contains the normal subgroup $image(H(q,F))$. Thus, if
we assume (A), we obtain
$$ n(\Gamma, F) = dim\, U(H(\Gamma, F)) $$
$$ \ge dim\, U(H(\Gamma', F)) + dim\, U(H(\Gamma'', F))$$
$$ = n(\Gamma', F) + n(\Gamma'', F)$$
from the previous lemma. But lemma \ref{lemma:s6} gives the reverse 
inequality, and this proves part (B). That (A) implies (C) is clear
by remark \ref{remark:3},  for we have
$$H(p,F)\circ i(\Gamma, F) = i(\Gamma'', F)$$
and
$$H(q, F) \circ i(\Gamma', F) = i(\Gamma, F) \circ q.$$

To check (A), replace $G$ in the lemma by the Zariski closure of 
$f(\Gamma)$, This makes $(G,j)$ an object of $C(\Gamma, F)$;
denote by $k: H(\Gamma, F) \to G_{min}$ the natural homomorphism.
Then $k \circ H(q, F)$ is the composite: $H(\Gamma',F)\to G\to G_{min}$.
By lemma \ref{lemma:a3}, this is an inclusion, This completes the
proof of the lemma.
\end{proof}

\begin{lemma}\label{lemma:s10} The $(G,j)$ in the previous lemma exists if 

(A) $\Gamma'' \cong \Z$, or

(B) $\Gamma''$ is a abelian torsion group.

\end{lemma}

\begin {proof} Case (A): Here $\Gamma$ is a semidirect product.
Choose $\gamma\in \Gamma$ that maps to a generator of $\Gamma''$.
Then $\sigma(\delta) = \gamma\delta \gamma^{-1}$, for $\delta\in \Gamma'$
gives an automorphism of $\Gamma'$, and induces therefore an 
automorphism $H(\sigma, F)$ of the algebraic group $H(\Gamma', F)$.
By lemma \ref{lemma:a4}, $A=Aut\,H(\Gamma, F)$
is an algebraic group. Let $G$ be the semidirect product of $A$ and 
$H(\Gamma, F)$, and define $j:\Gamma \to G(F)$ by $j(\delta) = 
i(\Gamma')(\delta)$ for $\delta \in \Gamma'$ and $j(\gamma) = 
H(\sigma,F)$.

Case (B):  First assume that $\Gamma'' $ is finite. Let 
$\rho:H(\Gamma', F) \to GL(V)$ be a faithful representation where $V$ 
is finite dimensional vector space defined over $F$. Consider the 
induced representation
$W = F[\Gamma]\otimes_{F[\Gamma']}V$. Let  $G = GL(W)$, and 
$j:\Gamma\to G$ be the action of $\Gamma$ on $W$. Clearly, there is a 
monomorphism $k: H(\Gamma', F) \to GL(W)$ so that $k\circ i(\Gamma', 
F) = j|_{\Gamma'}$, so the result follows.

In the general case, let 
$${\mathcal S} = \{\pi \subseteq \Gamma \, |\, \pi \supseteq \Gamma', \,
\pi/\Gamma'\> is \> finite\}$$
If $\pi_{1},\pi_{2}\in {\mathcal S}$ and $\pi_1 \subseteq \pi_2$, then 
$\pi_1/\pi_2$ is finite. Therefore, from the first part of the 
previous lemma, $H(\pi_1, F) \to H(\pi_2, F)$ has no kernel. From 
lemma \ref{lemma:a6}, we get:
$$\Gamma = \lim_{\stackrel{\longrightarrow}{\pi\in {\mathcal S}}}\, \pi \to
          \lim_{\stackrel{\longrightarrow}{\pi\in {\mathcal S}}}\, 
         H(\pi,F)(F) \to G(F)$$
and this completes the proof.
\end{proof}

We can now prove theorem \ref{thm:solv}:

\begin{proof} Part (A) has already proved in lemma \ref{lemma:s8}. 
 The theorem is certainly true if $\Gamma \cong \Z$, it is also true
 if $\Gamma$ is a commutative torsion group by lemma \ref{lemma:s7}.
For the general case,  we note that $\Gamma$ has a filtration:
$$\Gamma = \Gamma_0\supset \Gamma_1 \supset \ldots \supset 
\Gamma_{m+1} = \{1\}$$
such that the successive quotients are either abelian torsion 
groups or isomorphic to $\Z$. We proceed by induction on $m$, and so 
we can assume the theorem for $\Gamma_{1}$. Now lemmas \ref{lemma:s9},
\ref{lemma:s10} give the theorem for $\Gamma$.
\end{proof}

\begin{lemma} If $\Gamma$ is solvable of finite rank, then 
$H(\Gamma,F)\to H(\Gamma)_{F}$ an isomorphism.
\end{lemma}

\begin{proof} Put $G = H(\Gamma)_F$ and let $f :\Gamma \to 
H(\Gamma)_F(F)$ be given by $f(\gamma) = i(\Gamma)(\Gamma)$ for 
$\gamma\in \Gamma$. Then $(G,f)\in Obj\, C(\Gamma, F)$ and $dim\, 
U(G) = dim\, U(H(\Gamma)) = rk\Gamma$ by the theorem. From lemma
\ref{lemma:s8}, $H(\Gamma, F)\to G_{min} = G$ is an isomorphism,
and so the lemma is proved.
\end{proof}

\begin{thm}\label{thm:compltn}
Let $\Gamma$ be a finitely generated solvable group of finite rank.
There is a finite set of prime numbers $S$ so that 
$$H(\Gamma)_{\Q_l} \to H(\hat\Gamma, \Q_l)$$
is an isomorphism for all primes $l\notin S$.
\end{thm}

\begin{proof}
We may regard $H(\Gamma)$ as an algebraic subgroup of $(GL_n)_\Q$.
Because $\Gamma$ is finitely generated, $i(\Gamma)(\Gamma) \subset
GL_n(S^{-1}\Z)$ for some finite set of primes $S$. If $l$ is a prime 
not in $S$, then $GL_n(S^{-1}\Z) \subset GL_n(\Z_l)$, and the latter
is a profinite group. This gives a continuous homomorphism
$f_l:\hat \Gamma \to GL_n(\Z_l)$. Because $image(i(\Gamma))\subset
H(\Gamma)(\Q_l)$ and the second group  is closed  in $GL_n(\Q_l)$, 
we see that $f_l(\hat \Gamma) \subset H(\Gamma)(\Q_l)$.
The object $(H(\Gamma)_{\Q_l}, f_l)$ of $C(\hat \Gamma, \Q_l)$ gives 
an epimorphism
$H(\hat \Gamma, \Q_l) \to H(\Gamma)_{\Q_l}$. By the previous lemma
and remark \ref{remark:1}, we have $H(\Gamma)_{\Q_l} \to H(\hat \Gamma, 
\Q_l)$, and these arrows are inverses of each other.
\end{proof}

\begin{lemma}\label{lemma:s12}
If $\Gamma$ is solvable of finite rank and if $H(\Gamma)^o$ is 
unipotent, then there are normal subgroups $\Gamma_1\supseteq \Gamma_2$
of $\Gamma$ so that 

(a) $\Gamma/\Gamma_1$ is finite,

(b) $\Gamma_1/\Gamma_2$ is nilpotent, and

(c) $\Gamma_2$ is torsion.
\end{lemma}

\begin{proof} We take $\Gamma_1 = i(\Gamma)^{-1}H(\Gamma)^o(\Q)$
and $\Gamma_2 = ker(i(\Gamma))$. Then $\Gamma_1/\Gamma_2 \subset 
H(\Gamma)^o(\Q)$ and the latter is nilpotent. From theorem 
\ref{thm:solv}
$\Gamma_2$ is torsion, and $\Gamma/\Gamma_1 \subset 
H(\Gamma)/H(\Gamma)^o(\Q)$ is finite.
\end{proof}


\section{Fundamental groups of varieties}

Let $Y$ be a normal variety defined over a subfield $K\subset \C$
with a $K$-rational point $y_0$. For any field extension
$K'\supseteq K$, set $Y_{K'} = Y\times_{spec\,K} spec\,K'$.
Let $\pi_1^{alg}(X)$ denote the algebraic fundamental 
group of a connected scheme $X$ (with an unspecified base point),
 and $\hat \pi$ the profinite completion of a group
$\pi$.  Then we have a split exact sequence
$$1\to \hat \pi_1(Y_\C^{an})\to \pi_1^{alg}(Y_{K})\to Gal(\bar K/K)\to
1$$
(where the splitting depends on $y_0$)
\cite[IX 6.4, XII 5.2]{sga1}. Thus $Gal(\bar K/K)$ acts
continuously on $\hat \pi_1(Y_\C)$, and therefore also
on the abelianization  $H_1(Y_{\C}^{an},\Z)\otimes \hat\Z$
and its pro-$l$ part $ H_1(Y_\C^{an}, \Z)\otimes \Z_l$.

\begin{lemma}\label{lemma:finiteness}
Let $Y$ be a normal variety defined over a finitely generated field
$K\subset \C$. Then  $H_0(Gal(\bar K/K), H_1(Y_\C^{an}\otimes \Z_l))$ is finite.
\end{lemma}

\begin{proof} 
$H_1(Y)/(torsion)\, \otimes \Z_l$ is dual to 
$H_{et}^1(Y_{\bar K},\Z_l)$ as a $Gal(\bar 
K/K)$-module, thus it suffices to prove that the second group has
no invariants. Let $p:\tilde Y\to Y_{K'}$
be a desingularization defined over a finite extension $K'\supseteq K$.
By Zariski's main theorem the geometric fibers of $p$ are connected,
thus $H_{et}^1(Y_{\bar K},\Z_l)$ injects into
$H_{et}^1(\tilde Y_{\bar K},\Z_l)$, and
this is compatible with the $Gal(\bar K/K')$-action.
$H^1(\tilde Y_{\bar K})$ has no $Gal(\bar K/K')$-invariants, because the 
eigenvalues of the Frobenius at any prime of good reduction have
absolute value $q^{1/2}$ or $q$ 
by \cite[sect. 3.3]{Deligne}.
\end{proof}

\begin{remark} The argument can be simplified (and lengthened) in a 
couple of ways.
An appropriate Lefschetz type argument allows one to reduce to the case
where $Y$ is a curve where the relevant estimate on eigenvalues of the
Frobenius goes back to Weil.
Alternatively, when $Y$ is curve, one can deduce the finiteness of  
$H_0(Gal(\bar K/K), H_1(Y_\C^{an}\otimes \Z_l))$ directly from class
field theory.
\end{remark}

\begin{thm}
Let $X$ be a normal (not necessarily complete)  algebraic variety  defined
over $\C$. Let $\pi = D^0\pi \supseteq D^1\pi \supseteq \ldots$ be the
derived series of $\pi = \pi_1(X,x_0)$. If there is  a natural number
$n$ so that $\pi/D^n\pi$ is solvable of finite rank, then there are 
normal subgroups $P\supseteq Q \supseteq D^n\pi$ of $\pi$ so that 

(a) $\pi /P$ is finite,

(b) $P/Q$ is nilpotent, and

(c) $Q/D^n\pi$ is a torsion group.
\end{thm}

\begin{proof}
Put $\Gamma = \pi/D^n\pi$, and $T = H(\Gamma)^0/U(H(\Gamma))$.
By lemma \ref{lemma:s12}, the theorem follows once it has been proved 
that $T$ is trivial. We may assume that $X$ and $x_0$ are defined 
over a finitely generated field $K\subset \C$. There is an
 action of $Gal(\bar K/K) $ on $\hat \pi$, and also on
$\hat \Gamma$, because this is a quotient of $\hat \pi$ by the closure
of $D^n(\hat \pi)$.
Choose a prime $l$, so that $\Gamma\to H(\Gamma)(\Q)$ extends to a 
continuous homomorphism $\hat \Gamma \to H(\Gamma)(\Q_l)$. For
such a prime, $H(\hat \Gamma, \Q_l) = H(\Gamma)_{\Q_l}$ by
theorem \ref{thm:compltn}. Thus the Galois action on $\hat \Gamma$
yields a homomorphism from $Gal(\bar K/K)$ to the group of 
$\Q_l$-rational points of $G = Aut(H(\Gamma))$. After replacing $K$
by a finite extension, if necessary, we can assume 
$image(\rho) \subset G^0(\Q_l)$. By lemma \ref{lemma:a5}, the action
of $Gal(\bar K/K)$ on $T_{Q_l}$ is trivial. Let 
$\pi' = ker[\pi \to (H(\Gamma)/H(\Gamma)^0)(\Q)]$
Then $\pi' = \pi_1(Y,y_0)$ where $Y$ is an etale cover of $X$.
The composite
$$\hat \pi' \hookrightarrow  \hat \pi \to 
(H(\Gamma)/U(H(\Gamma))(\Q_l)$$
factors through
$$H_1(Y,\Z)\otimes \hat \Z = \hat \pi_{ab}' \to T(\Q_l),$$
$T(\Q_l)$ contains an open pro-$l$-group, thus $\hat \pi'$
further factors through:
$$h:H_1(Y)\otimes\hat \Z \to H_1(Y)\otimes \Z_l \oplus A \to T(\Q_l)$$
where $A$ is a finite group. However  
$H_0(Gal(\bar K/K), H_1(Y)\otimes \Z_l)$  is finite by lemma
\ref{lemma:finiteness}.
So we deduce that the image of $H_1(Y)\otimes \Z_l$ in $T(\Q_l)$
is finite, because $Gal(\bar K/K)$ acts trivially on $T(\Q_l)$.
However the image $h$ is Zariski dense. Thus $T$ is finite and
connected, and therefore trivial. This proves the theorem.

\end{proof}

\begin{cor} If the fundamental group of a normal complex variety 
is solvable and possesses a faithful representation into $GL_n(\Q)$,
then it is virtually nilpotent i.e. it must contain a nilpotent
subgroup of finite index.
\end{cor}

\begin{proof} We can assume that the  fundamental group is torsion free
after passing to a subgroup of finite index \cite[lemma 8]{S}.
The theorem implies that this must contain a nilpotent group
of finite index.
\end{proof}

\section{Fundamental groups of K\"ahler manifolds}

A group $\Gamma$ will be called quasi-K\"ahler if it there
exists a connected compact K\"ahler manifold $X$, and a divisor
with normal crossings $D\subset X$, such that $\Gamma\cong \pi_1(X-D)$.
(Note that by resolution of singularities \cite{AHV, BM},
 it is enough to assume that $D$ is an analytic subset.) 

The proof of the following lemma will be given in the appendix.

\begin{lemma} A subgroup of  a {\qk} group of finite index is {\qk}.
\end{lemma}

\begin{lemma} Let $A$ be a finitely generated abelian group, and $M$
a nontrivial one dimensional $\C[A]$-module. Then $H^i(A,M)=0$ for all
$i$.
\end{lemma}

\begin{proof} This is clear for cyclic groups by direct computation.
In general, express $A$ as a product of cyclic groups $\prod_{i}{C_i}$ and 
$M$ as a tensor product of $C_i$-modules, and apply the K\"unneth 
formula.
\end{proof}

Set $\Gamma^{ab} = \Gamma/D\Gamma$.

\begin{lemma}\label{lemma:h1} Suppose that $\Gamma$ is a finitely
generated group
and $A=\Gamma/N$ an abelian quotient.
Suppose that $M$ is a nontrivial one dimensional $A$-module then
$$H^1(\Gamma,M) \cong Hom_{\Z[A]}(N^{ab}, M)\cong 
Hom_{\Q[A]}(N^{ab}\otimes\Q, M).$$
\end{lemma}

\begin{proof} From the Hochschild-Serre spectral sequence associated
to the extension $1\to N \to \Gamma \to A\to 0$, we obtain an exact sequence
$$0\to H^1(A, M)\to H^1(\Gamma, M)\to H^0(A,H^1(N,M))\to H^2(A,M).$$ 
By the previous lemma, this gives an isomorphism
$$H^1(\Gamma,M)\cong H^0(A,H^1(N,M))$$
Furthermore 
$H^1(N,M) \cong Hom(N,M) = Hom(N^{ab},M)$ as $A$-modules. Therefore
$H^0(A,H^1(N,M)) \cong Hom_{\Z[A]}(N^{ab}, M)$.
\end{proof}

Given a character  $\rho\in Hom(\Gamma, \C^*)$, let $\C_\rho$ denote
the associated $\Gamma$-module.
We define $\Sigma^1(\Gamma) $ to be the set of
characters $\rho\in Hom(\Gamma,\C^*)$ such  that $H^1(\Gamma,\C_\rho)$ 
is nonzero. Let us say that a $\Gamma$-module  $V$ is quasiunipotent if
there is a subgroup $\Gamma'\subseteq \Gamma$ of finite index whose 
elements act unipotently  on $V$.

\begin{lemma} Let $A$ be a finitely generated abelian group and $V$
a finite dimensional $\C[A]$-module. Then $A$ acts quasiunipotently on
$V$ if and if the only characters $\rho\in Hom(A,\C^*)$ for which
$Hom_{\C[A]}(V,\C_\rho)\not= 0$ are torsion characters (i.e. elements of
$Hom(\Gamma, \C^*)$ of finite order).
\end{lemma}

\begin{proof} 
Define $V_\rho$ to be the generalized eigenspace associated to
a character $\rho$. In other words, $V_\rho$ is the maximal subspace
on which $a-\rho(a)$ is nilpotent for all $a\in A$. $V$ is a direct 
sum of these eigenspaces, thus we can assume that $V=V_\rho\not= 0$.
To complete the proof, observe that $Hom(V_\rho,\C_{\rho'})\not= 0$
if and only if $\rho=\rho'$, and that $V_\rho$ is quasiunipotent if
and if $\rho$ is torsion.
\end{proof}

\begin{lemma}\label{lemma:torsion} Let  $\Gamma$ be a finitely 
generated group and $\Gamma'\subseteq \Gamma$ a subgroup of finite
index such that $V= (\Gamma'\cap D\Gamma)^{ab}\otimes \Q$ is finite 
dimensional. Then $\Gamma'$
acts quasiunipotently on $V$ if and only if
 $$\Sigma^1(\Gamma')\cap image(Hom(\Gamma,\C^*)\to Hom(\Gamma',\C^*))$$
consists of torsion characters.
\end{lemma}

\begin{proof} 
Set $N = \Gamma'\cap D\Gamma$. Then $\Gamma'/N$ is
isomorphic to the image of $\Gamma'$ in $\Gamma^{ab}$, and therefore
$Hom(\Gamma'/N,\C^*)$ is coincides with
$image(Hom(\Gamma,\C^*)\to Hom(\Gamma',\C^*)).$
Thus, lemma \ref{lemma:h1} implies that 
$$S=(\Sigma^1(\Gamma')-\{1\})\cap image(Hom(\Gamma,\C^*)\to Hom(\Gamma',\C^*))$$
is the set of nontrivial characters $\rho\in Hom(\Gamma'/N,\C^*)$ 
for which $Hom_{\C[A]}(V\otimes\C,\C_\rho)\not= 0$.
Thus $\Gamma'/N$ acts quasiunipotently on $V$ if and only if $S$ 
consists of torsion characters by the previous lemma.
The $\Gamma'$ action on $V$ factors through $\Gamma'/N$, thus
the lemma is proved.
\end{proof}

\begin{lemma}\label{lemma:kron} Let $K\supseteq \Q$ be a finite  extension,
and $O_K$ the ring of integers of $K$. Let $\{\sigma_1,\ldots\sigma_n\}$ 
be the set of all embeddings $K$ into $\C$.
 Then for any finitely generated group $\Gamma$,
  $$Hom(\Gamma, U(1))\cap \bigcap_{i}\sigma_i\circ Hom(\Gamma, O_K^*) $$ 
  consists of torsion characters.
\end{lemma}

\begin{proof} This follows from Kronecker's theorem that an algebraic 
integer is a root of unity if and only if all its Galois conjugates 
have absolute value one.

\end{proof}

\begin{thm}\label{thm:qunip}  Let  $\Gamma$ be a {\qk} group such that 
$D\Gamma$  is a finitely generated. Then for any subgroup 
$\Gamma'\subseteq \Gamma$ of finite index, $\Gamma'$ acts 
quasiunipotently on the finite dimensional vector space
$(\Gamma'\cap D\Gamma)^{ab}\otimes \Q$.
\end{thm}

\begin{proof} By  lemma \ref{lemma:torsion}, it is enough to show 
that the intersection $S$ of  
$\Sigma^1(\Gamma')$ and the image of $Hom(\Gamma,\C^*)$
 consists of torsion characters.
 The subgroup $\Gamma'\cap D\Gamma\subseteq D\Gamma$
 has finite index, and is therefore  finitely generated. Thus
the set of characters of $\Gamma'$
which correspond to one dimensional quotients
of $(\Gamma'\cap D\Gamma)^{ab}\otimes \C$ are defined over  the ring of 
integers of a finite extension $K$ of $\Q$. It then follows from lemma 
\ref{lemma:h1} that $S$ is a finite subset of 
$Hom(\Gamma, O_K^*)$. $S$ is evidently stable under $Aut(\C)$ and thus
lies in $\cap_{\sigma:K\hookrightarrow\C}\, \sigma\circ Hom(\Gamma,O_K^*)$.

Theorem V.1.6 of \cite{A2} implies that 
$\Sigma^1(\Gamma')$ is a finite union of translates of subtori
of $Hom(\Gamma')$ by unitary characters. $S$, which is the 
intersection of this set with $image(Hom(\Gamma,\C^*))$, must clearly 
inherit a similar structure. In particular, being finite, it follows 
that $S$ consists of unitary characters.
Therefore the theorem follows from lemma \ref{lemma:kron}.

\end{proof}

\begin{lemma}\label{lemma:key} Let  $\Gamma$ be a solvable group of 
finite rank. 
Suppose that  every subgroup $\Gamma'\subseteq \Gamma$ of finite index
acts quasiunipotently on the finite dimensional vector space
$(\Gamma'\cap D\Gamma)^{ab}\otimes \Q$. Then there exists
 normal subgroups $\Gamma_1\supset \Gamma_2$
of $\Gamma$ so that 

(a) $\Gamma_1$ has finite index,

(b) $\Gamma_1/\Gamma_2$ is nilpotent, and

(c) $\Gamma_2$ is torsion.
\end{lemma}

\begin{proof}  It suffices to prove that  $G= H(\Gamma)^o$ is
unipotent by  lemma \ref{lemma:s12} (then $\Gamma_1$ can be be 
taken to be $i(\Gamma)^{-1}(G)$  and $\Gamma_2= ker\, i(\Gamma)$).
The unipotency of $G$ will follow from
lemma \ref{lemma:unip}, once we show that the action
of $G/DG$ on $DG/D^2G$, by conjugation,  
is unipotent. 
The map 
$$(\Gamma_1\cap D\Gamma)^{ab}\otimes\Q \to DG(\Q)/D^2G(\Q)$$ 
is compatible 
with the $\Gamma_1$-actions, and is surjective, because
the image of $\Gamma_1\cap D\Gamma$ is Zariski dense in  $DG(\Q)$.
By hypothesis, $\Gamma_1$ contains a finite index subgroup
$\Gamma''$ which acts unipotently on $(\Gamma_1\cap D\Gamma)^{ab}\otimes\Q$.
As $\Gamma''$ is Zariski dense in $G$, the lemma follows.
\end{proof}

\begin{thm}
 Let $\pi $ be a {\qk} group. Suppose that  $D\pi$ is  
 finitely generated and $\pi/D^n\pi$ is solvable of finite rank 
  for some natural number $n$. Then there are 
normal subgroups $P\supseteq Q \supseteq D^n\pi$ of $\pi$ so that 

(a) $\pi /P$ is finite,

(b) $P/Q$ is nilpotent, and

(c) $Q/D^n\pi$ is a torsion group.
\end{thm}

\begin{proof} The  theorem  follows from
theorem  \ref{thm:qunip} and   lemma \ref{lemma:key}.

\end{proof}

\begin{cor} A  polycyclic {\qk} group is virtually nilpotent.
\end{cor}

\begin{proof} By \cite[4.6]{R}, a polycyclic group contains a
torsion free subgroup $\pi$ of finite index. The theorem implies
that $\pi$ must contain a  nilpotent subgroup of finite index.
\end{proof}

\appendix\section{Appendix: Construction of K\"ahler metrics}

Most of the  results in this appendix are  well known, but we 
indicate proofs for lack of a suitable reference. The following is elementary,
and left to the reader.

\begin{lemma}\label{lemma:posdef} Let $W= U\oplus V$ be a finite dimensional $\C$-vector 
space. Let  $Q_{i,t}$, $i=1,2$ be two hermitian forms  on
$W$ which depend continuously on a parameter $t$ varying over a compact
set $T$. Suppose that $Q_{1,t}$ is positive semidefinite and $Q_{1,t}|_{U}$ 
and $Q_{2,t}|_{V}$ are  positive definite for all $t$. Then there
exists a constant $x > 0$ such that $xQ_{1,t}+ Q_{2,t}$ is positive 
definite for all $t$.
\end{lemma}

\begin{lemma}\label{lemma:fs}
 Let $V$ be a holomorphic vector bundle over a complex manifold $X$. 
 Then the hyperplane
bundle $O(1)$ on $\pi:\PP(V)\to X$ carries a hermitian metric which
restricts to the Fubini-Study metric on every fiber.
This metric will be referred to as a Fubini-Study type metric.
\end{lemma}

\begin{proof} $O(1)$ is a quotient of $\pi^*V$.
Let $h$ be a hermitian metric on $V$. Then the metric on $O(1)$ 
induced from $\pi^*h$ has the desired properties.
\end{proof}

Given a hermitian metric $h$ on a line bundle, let $\tilde c_1(h)$ denote
the first Chern form, given locally by 
$(i/2\pi)\partial\bar\partial\log ||s||^{2}_{h}$ where $s$ is 
holomorphic section. 

\begin{lemma} Let $X$ be compact complex manifold with a K\"ahler form 
$\omega$. If $V$ is a holomorphic vector bundle and $h$ a
Fubini-Study type metric on $O_{\PP(V)}(1)$, then  
$C\pi^*\omega + \tilde c_1(h)$ is a K\"ahler
form on $\pi:\PP(V)\to X$ for all $C >> 0$.
\end{lemma}

\begin{proof}  $C\pi^*\omega + \tilde c_1(h)$ is a 
real $(1,1)$-form, and it is positive if $C >>0$ by \ref{lemma:posdef}.
\end{proof}

\begin{lemma} The blow up of a compact K\"ahler manifold along a 
submanifold is K\"ahler.
\end{lemma}

\begin{proof} 
Let $\pi:\tilde X \to X$ be the blowup of $X$ along a closed 
submanifold $S$, and let $E$ be the exceptional divisor.
By construction \cite{gh},
there is an open tubular neighbourhood 
$U$ of $S$, such that the preimage $\tilde U$ embeds into
$U\times \PP(N)$, where $N$ is the normal bundle of $S$.
The restriction of $O(-E)$ to $U$ coincides with the pullback of 
$O(1)$. Let $h_1$ be the restriction of a Fubini-Study type metric
to  $O(-E)|_{\tilde U}$. Choose a $C^{\infty}$ cutoff function $\rho :\tilde 
X\to [0,1]$ which vanishes outside $\tilde U$
and is identically $1$ on a neighbourhood of $V$ of $E$.
Let $h_2$ be a constant metric on the trivial bundle $O(-E)|_{X-\bar 
V}$. Then $h= \rho h_1 + (1-\rho)h_2$ defines a  metric on $O(-E)$.
If $\omega$ is a K\"ahler metric on $X$, then lemma \ref{lemma:posdef} 
shows that $C\pi^{*}\omega + \tilde c_1(h)$ is K\"ahler for $C >> 0$.

\end{proof}

\begin{prop} Suppose that $X$ is a compact K\"ahler manifold,
  $D\subset X$ is a divisor with normal crossings, and $Y^o\to X-D$
is a finite sheeted unramified cover. Then $Y^o$ has a (nonsingular)
 K\"ahler compactification $Y$, such that $Y- Y^o$ is a divisor
with normal crossings.
\end{prop}

\begin{proof} Let $X^o = X-D$ and
let $n$ be the number of sheets of $Y^o\to X^o$. 
Then there exists a principle $S_n$-bundle $P$, such that
 $Y^o$ is isomorphic to $ P\times_{S_n} \{1\ldots n\}$
 Let $S_n\to GL_n(\C)$ be the permutation
representation associated to the standard basis $\{e_i\}$. Then we 
obtain a holomorphic vector bundle $V^o = P\times_{S_n} \C^n$.
There is a holomorphic embedding $Y^o\mapsto V^o$ induced by
the map $\{1\ldots n\} \to \C^n$ given by $i\mapsto e_i$.
$V^o$ is a flat vector bundle, so it has a natural flat connection $\nabla$. 
$V^o$   extends to a holomorphic bundle
$V$ on $X$ with regular singularities with respect to $\nabla$
\cite{RegSing}.  Let $Y'$ be the closure  of $Y^o$ in $V$.
It is easy to check that $Y'$ is an analytic subset of $V$.
We can embed $V$ into the projective space bundle $P = \PP(V\oplus O)$.
By resolution of singularities \cite{AHV, BM},
 there exists a commutative diagram
$$\begin{array}{ccc}
Y & \hookrightarrow & \tilde P \\
\downarrow&& \downarrow \\
Y' & \hookrightarrow & P\\
\end{array}
$$
where $\tilde P \to P$ is a composite of blow ups along smooth 
centers lying over $D$, $Y$ is nonsingular, and $Y-Y^o$ is a divisor
with normal crossings. $\tilde P$ is K\"ahler by the previous lemmas, 
therefore the same is true of $Y$.
 
\end{proof}

As an immediate corollary, we obtain:

\begin{lemma} A  finite index subgroup of a {\qk} group is {\qk}.
\end{lemma}


\end{document}